\documentstyle[12pt]{amsart} 
\let\cal\mathcal

\makeatletter
\def\displaylines#1{\displ@y
  \halign{\hbox to\displaywidth{$\@lign\hfil\displaystyle##\hfil$}\crcr
    #1\crcr}}
\makeatother

\numberwithin{equation}{section}
\newcommand{\beq}{\begin{equation}}
\newcommand{\eeq}{\end{equation}}

\newtheorem{theorem}{Theorem}[section]

\newtheorem{prop}[theorem]{Theorem}

\newtheorem{coroll}[theorem]{Corollary}

\newtheorem{lemm}[theorem]{Lemma}

\newtheorem{anoprop}[theorem]{Proposition}

\newtheorem{eg}{\rm\sl \uppercase{Example}}[section]

\begin{document}

\title{Excision in Banach simplicial and cyclic cohomology}
\author{Zinaida A. Lykova}
\address{Zinaida A. Lykova\\
 Department of Mathematics and Statistics\\
 Fylde College, Lancaster University \\
 Lancaster LA1 4YF \\
 England}
\email{maa032@@cent1.lancs.ac.uk}

\begin{abstract}

   We prove that, for every extension of Banach algebras
$  0 \rightarrow B \rightarrow  A  \rightarrow D \rightarrow  0 $  such that 
 $B$ has a left or right bounded approximate identity,
 the existence of an associated long exact sequence of Banach simplicial or cyclic
 cohomology groups is equivalent to the  existence of one for homology groups.
It follows from the continuous version of a result of Wodzicki that 
associated long exact sequences exist.
 In particular, they exist for every  extension of 
$C^*$-algebras.

\end{abstract}


\thanks{I am indebted to  the Mathematics Department of 
the University of California and the Mathematical Sciences Research 
Institute at Berkeley  for hospitality while this work was carried out.
 Research at MSRI is supported in part by NSF grant DMS-9022140.}

\maketitle

\section{Introduction}

  Significant results of A. Connes on non-commutative differential geometry
[Co1] have led  much research interest to  the computation of cyclic
(co)homology groups in recent years; see, for example, [Lo] and  [Co2] for many 
references and [ChS], [Gr], [He2], [Wo1] and [Wo2] for the continuous theory 
of these groups. A promising approach to the calculation of cyclic cohomology
groups is to break it down by making use of extensions of Banach algebras; 
this is a standard device in the study of  various
properties of Banach algebras.  Here we establish (Theorem 3.5)
that, for every weakly admissible extension of Banach algebras
$  0 \rightarrow B \rightarrow  A  \rightarrow D \rightarrow  0 $,
 the existence of an associated long exact sequence of Banach simplicial 
 cohomology groups is equivalent to the  existence of one for homology groups.
Theorem 3.8 shows that the same result is true for cyclic (co)homology groups when $B$
has a left or right bounded approximate identity.
Recall that the existence of long exact sequences of simplicial and cyclic
{\em homology} groups  associated with an extension of  algebras 
(as distinct from cohomology and,  mainly, in an algebraic context) 
has been studied by M. Wodzicki [Wo2], [Wo3]; see also [Lo].  
Wodzicki remarked that his result could easily be extended to the continuous 
case under some hypotheses on the extension 
(see Remark (3) and Corollary 4 
[Wo2] and Remark 8.5(2) [Wo3]).

    From Theorems 3.5, 3.8 and the continuous version of 
Wodzicki's result  (Theorem 4.1)
for homology groups we deduce the existence of long exact 
sequences of  simplicial and cyclic
cohomology groups 
for every extension when $B$ has a left or right 
bounded approximate identity. In particular, this is true for every  extension of 
$C^*$-algebras. This gives 
an effective tool for computing Banach simplicial
and cyclic cohomology groups (see Proposition 4.4). 
We apply it to some natural classes of Banach algebras.

  Finally, I am grateful to M.Wodzicki for suggestions about the continuous homology of
Banach algebras.

\section{Definitions and notation}

  We recall some notation and terminology used in the homological theory of
Banach algebras.
 
    Throughout the paper $id$ denotes the identity operator. We denote the
projective tensor product of Banach spaces by  $\hat{\otimes}$ 
(see, for example, [Pa]). Note that by 
$Z^{\hat{\otimes} 0} \hat{\otimes} Y$ we mean $Y$ and by 
$Z^{\hat{\otimes} 1}$ we mean $Z$.

 Let $A$ be a Banach algebra, not necessarily unital.
 We will define the Banach version of the cyclic homology ${\cal HC}_n(A)$ 
(see, for example, [Lo], [Wo2] or [He2]).
     We denote by  $ C_n(A),~ n = 0, 1,...,$ the $(n+1)$ fold projective tensor power 
$A^{ \hat{\otimes}(n+1)} = A \hat{\otimes} \dots \hat{\otimes} A$ of $A$;
 we shall call the elements of this Banach space 
{\it $n$-dimensional chains}. We let
$ t_{\underline{n}}: C_n(A) \rightarrow  C_n(A), \; n = 0, 1,...,$
denote the operator given by
$t_{\underline{n}} (a_0 \otimes a_1\otimes \dots \otimes a_n) = 
(-1)^n (a_n \otimes a_0\otimes \dots \otimes a_{n-1}),$ 
and we set $t_{\underline{0}} = id$. We let $ CC_n(A)$ denote the 
quotient space of $C_n(A)$ modulo the closure of the linear span of elements of the form 
$ x -  t_{\underline{n}}x$ where $ n = 0, 1, \dots .$ Note that, by
Proposition 4 [He2],
${\rm Im}~ (id_{C_n(A)} - t_{\underline{n}})$ is closed in $C_n(A)$ and
so $ CC_n(A) = C_n(A)/{\rm Im}~ (id_{C_n(A)} - t_{\underline{n}})$. 
We also set 
$ CC_0(A) = C_0(A) = A.$

     From the chains we form the standard homology complex 

\hspace{2.5cm}
$ 0 \leftarrow C_0(A) \stackrel {d_0} {\leftarrow} \dots 
\leftarrow C_n(A) 
 \stackrel {d_n} {\leftarrow}  C_{n+1}(A) \leftarrow \dots , 
\hfill {({\cal C}_{\sim} (A))}$
\vspace*{0.2cm}
\newline
where the differential $d_n$ is given by the formula
\[
d_n (a_0\otimes a_1\otimes \dots \otimes a_{n+1}) = 
\]
\[
 \sum_{i=0}^n (-1)^i (a_0 \otimes \dots \otimes a_i a_{i+1} \otimes \dots \otimes a_{n+1}) +
(-1)^{n+1} (a_{n+1} a_0 \otimes \dots \otimes a_n).
\]
It is not difficult
to verify that these $d_n$ induce  operators
$dc_n: CC_{n+1}(A) \rightarrow CC_n (A) $ 
in the respective quotient spaces. Thus we obtain a quotient complex 
 ${\cal CC}_{\sim} (A)$ of the complex ${\cal C}_{\sim} (A)$.
 The $n$-dimensional homology of ${\cal C}_{\sim} (A)$, denoted
by  ${\cal H}_n(A)$, is called the {\it $n$-dimensional Banach  
simplicial homology group} of the Banach algebra $A$. The $n$-dimensional 
 homology of ${\cal CC}_{\sim} (A)$, denoted
by   ${\cal HC}_n(A)$, is called the {\it $n$-dimensional Banach cyclic 
homology group} of $A$.

   We  also consider the  complex 

\hspace{2.5cm}
$ 0 \stackrel {d r_{-1}}{ \leftarrow} C_0(A) \stackrel {d r_0} 
{\leftarrow} \dots \leftarrow C_n(A) 
 \stackrel {d r_n} {\leftarrow}  C_{n+1}(A) \leftarrow \dots , 
\hfill {( {\cal CR}_{\sim} (A))}$
\vspace*{0.2cm}
\newline
where the differential $d r_n$ is given by the formula
\[
d r_n (a_0\otimes a_1\otimes \dots \otimes a_{n+1}) = 
 \sum_{i=0}^n (-1)^i (a_0 \otimes \dots \otimes a_i a_{i+1} 
\otimes \dots \otimes a_{n+1})
\]
and we set $d r_{-1} (a_0) = 0.$  The $n$-dimensional 
 homology of ${\cal CR}_{\sim} (A)$, denoted
by   ${\cal HR}_n(A)$, is called the {\it $n$-dimensional Banach Bar 
homology group} of $A$.

   Note that, by definition, $dc_0 = d_0$, so that  ${\cal HC}_0(A) = 
{\cal H}_0(A) = A/{\rm Im}~ d_0$.

  Obviously   ${\cal H}_n(A)$ is just 
another way of writing  the Hochschild homology group ${\cal H}_n(A, A)$
(see [He1; II.5.5]).

  For   a Banach $A$-bimodule $X$,
 we will denote  the $n$-dimensional Banach
cohomology group of $A$ with coefficients in $X$ by $ {\cal H}^n(A,X)$
(see, for example, [Jo1] or [He1]).
  Recall that a Banach $A$-bimodule $M = (M_*)^*$, where $M_*$ is a Banach
$A$-bimodule, is called {\it dual}. Here  $E^*$ is the dual Banach space of a Banach
space  $E$. A Banach algebra $A$ such that 
 $ {\cal H}^1(A,M) =\{ 0 \}$
for all dual $A$-bimodules $M$ is called {\it amenable}.

 The $n$th cohomology of the dual complex ${\cal C}^{\sim}(A)
 \stackrel {def} {=}  
{\cal C}_{\sim} (A) ^*$, denoted
by  ${\cal H}^n(A)$, is called the {\it $n$-dimensional Banach  
simplicial cohomology group} of the Banach algebra $A$. The $n$th
 cohomology of the dual complex  ${\cal CC}^{\sim}(A) 
\stackrel {def} {=}
{\cal CC}_{\sim} (A) ^* $, denoted
by   ${\cal HC}^n(A)$, is called the {\it $n$-dimensional Banach cyclic 
cohomology group} of $A$  (see [Co1] or [He2]).
 The $n$th cohomology of the dual complex
 ${\cal CR}^{\sim}(A) \stackrel {def} {=}
{\cal CR}_{\sim} (A) ^*$, denoted
by  ${\cal HR}^n(A)$, is called the {\it $n$-dimensional Banach  
Bar cohomology group} of the Banach algebra $A$.

    Note that, by definition,  ${\cal HC}^0(A) = 
{\cal H}^0(A)$ coincides with the space $A^{tr} =
\{ f \in A^* : f(ab) = f(ba) \; {\rm for \; all} \; a, b \in A \}$ of continuous
traces on $A$.

 The canonical identification of $(n +1)-$linear functionals on $A$
 and $n-$linear operators from $A$ to $A^*$ shows that   ${\cal H}^n(A)$ is just 
another way of writing   ${\cal H}^n(A, A^*)$.

  The vanishing of  ${\cal HR}^n(A)$ for a Banach algebra $A$ implies
the existence of the Connes-Tsygan exact sequence for $A$. One can find the conditions
for this in Theorems 15 and 16 [He2].

 The short exact sequence of Banach 
spaces and  continuous linear operators
\[
 0 \rightarrow Y \stackrel {i} {\rightarrow}    Z  \stackrel {j} {\rightarrow} W
 \rightarrow  0
\]
is called {\it admissible} if there exists a continuous operator 
$\alpha : W \rightarrow Z$ such that $ j \circ \alpha = id_W$.
   Recall that admissibility is equivalent to  the existence of  
 a continuous operator 
\newline
$\beta: Z \rightarrow Y$ such that $ \beta \circ i = id_Y$. This is also
 equivalent to  the existence of   continuous operators 
$\beta: Z \rightarrow Y$  and 
$\alpha : W \rightarrow Z$ such that  $ \beta \circ i = id_Y$, $ j \circ \alpha = id_W$ and
$  i \circ \beta  + \alpha \circ j = id_Z$.

 The short exact sequence of Banach 
spaces and  continuous linear operators
\newline
$ 0 \rightarrow Y \stackrel {i} {\rightarrow}    Z  \stackrel {j} {\rightarrow} W
 \rightarrow  0 $ 
is called {\it weakly admissible} if the dual short sequence 
\newline
$ 0 \rightarrow W^* \stackrel {j^*} {\rightarrow}    Z^*  \stackrel {i^*} {\rightarrow} Y^*
 \rightarrow  0$
is admissible.

 Definitions of a (co)chain complex, a morphism of complexes, the homology groups can be
found in any text book on homological algebra, for instance MacLane [Ma], Helemskii [He1]
(continuous case), Loday [Lo].

 A {\it  chain complex} ${\cal K}_{\sim} $ in the category of 
[Banach] linear spaces is a 
 sequence of [Banach] linear spaces and  [continuous] linear operators
\[
 \dots \leftarrow K_n \stackrel {d_n} {\leftarrow}  K_{n+1} 
\stackrel {d_{n+1}} {\leftarrow} K_{n+2}
 \leftarrow  \dots
\]
such that $d_n \circ d_{n+1} = 0 $ for every $n$.
 The {\it cycles} are the elements of $Z_n = {\rm Ker}~( d_{n-1}:
K_n \rightarrow  K_{n-1}).$ The   {\it boundaries} are the elements of 
$B_n = {\rm Im}~( d_n: K_{n +1} \rightarrow  K_n).$ The relation $d_{n-1} \circ d_n = 0 $
implies $B_n \subset Z_n.$ The {\it homology groups} are defined by $ H_n({\cal K}_{\sim}) = 
Z_n/B_n.$

  A {\it [continuous] morphism of chain complexes }
 ${\psi}_{\sim} : {\cal K}_{\sim} \rightarrow {\cal P}_{\sim}$ in the category of 
[Banach] linear spaces is a collection of  [continuous] linear operators 
${\psi}_n :  K_n \rightarrow P_n$ such that the following diagram is commutative for any $n$
\[
\begin{array}{ccc}
K_n & \stackrel {d_{n-1}} {\rightarrow}  &   K_{n-1}
\\
 \downarrow {\psi}_n & ~ &   \downarrow {\psi}_{n-1} 
\\
P_n & \stackrel {d'_{n-1}} {\rightarrow}  &   P_{n-1} .
\end{array}
\]
Such a morphism obviously induces a map   $H_n({\psi}_{\sim}) : H_n({\cal K}_{\sim} )
 \rightarrow H_n({\cal P}_{\sim})$.

  A [continuous] morphism of chain complexes 
 ${\psi}_{\sim} : {\cal K}_{\sim} \rightarrow {\cal P}_{\sim}$ of [Banach] linear spaces is {\it a
[topological] quasi-isomorphism} if  $H_n({\psi}_{\sim}) : H_n({\cal K}_{\sim} )
 \rightarrow H_n({\cal P}_{\sim})$ is  a
[topological] isomorphism for every $n $.

\section{ Connections between long exact sequences of homology and 
cohomology groups }

  To start with we need the following elementary result from the
homological algebra. 

{\sc 3.1 Lemma.}~~
 {\it Let 
\[
 \dots \rightarrow {\cal L}_{n-2}  \stackrel {\zeta_{n-1}} {\rightarrow} 
{\cal L}_{n-1}  \stackrel {\zeta_n} {\rightarrow}  {\cal L}_n 
 \stackrel {\zeta_{n+1}} {\rightarrow}  {\cal L}_{n+1} 
\stackrel {\zeta_{n+2}} {\rightarrow}  {\cal L}_{n+2}   \dots ,
\]
be an exact sequence of linear spaces and linear operators.
 Suppose that $\zeta_{n-1}$ and  $\zeta_{n+2}$ are isomorphisms.  
 Then $ {\cal L}_n = \{0\}$. }

 {\it Proof.} By the assumption that $\zeta_{n-1}$ is an isomorphism, so 
$ {\rm Im}~ \zeta_{n-1} = {\cal L}_{n-1} = {\rm Ker}~ \zeta_{n}.$
 Thus $ {\rm Im}~ \zeta_{n} = \{0\} $ and, by virtue of the exactness of the sequence,
$ {\rm Ker}~ \zeta_{n+1} =  {\rm Im}~ \zeta_{n} = \{0\}$.  By the assumption that $\zeta_{n+2}$ 
is an isomorphism, $ {\rm Ker}~ \zeta_{n+2} = \{0\} =  {\rm Im}~ \zeta_{n+1}$. Therefore, since 
$ {\rm Ker}~ \zeta_{n+1} = \{0\}$ and $ {\rm Im}~ \zeta_{n+1} = \{0\}$, the result follows. 
\qed

{\sc 3.2 Proposition.}~~ {\it Let  ${\cal K}_{\sim}$ and ${\cal P}_{\sim}$  be chain 
complexes of Banach
spaces and continuous linear operators. Suppose there is a continuous morphism  
of chain complexes ${\psi}_{\sim} : {\cal K}_{\sim} \rightarrow {\cal P}_{\sim}$ such that, for 
 each $n$, $\psi_n$ is injective and
${\rm Im}\psi_n$ is closed. Then 
 
   {\rm (I)} for  fixed $n$, if $H^k({\psi}_{\sim}^*)$ is an isomorphism
whenever $n +2 \ge k \ge n-1$, then $H_n({\psi}_{\sim})$ is a topological 
isomorphism.

 {\rm (II)} for  fixed $n$, if  $H_k({\psi}_{\sim})$ 
 is an isomorphism
whenever  $ n + 1 \ge k \ge n-2$, then  $H^n({\psi}_{\sim}^*)$  is a 
topological isomorphism.

   {\rm (III)} the following conditions are equivalent:

  {\rm (i)} ${\psi}_{\sim} $ is a topological quasi-isomorphism  of  the chain complexes
$ {\cal K}_{\sim} $ and ${\cal P}_{\sim}$; 

 {\rm (ii)} ${\psi}_{\sim}^* $ is a topological quasi-isomorphism  of the cochain, dual complexes
$ {\cal K}_{\sim}^* $ and ${\cal P}_{\sim}^*$. }

 {\it Proof.} By the assumption, there is a short exact sequence of complexes
\[
  0 \rightarrow   {\cal K}_{\sim}  \stackrel {{\psi}_{\sim}} {\rightarrow}     
 {\cal P}_{\sim}   \stackrel {{\phi}_{\sim}} {\rightarrow}    {\cal L}_{\sim}
 \rightarrow   0
\]
in the category of Banach spaces and continuous linear operators, where
$  {\cal L}_{\sim} $ is the chain  complex ${\cal P}_{\sim}/{\rm Im} {\psi}_{\sim}$
and $ \phi_n: P_n \rightarrow P_n/{\rm Im} {\psi}_{n}$ is the quotient mapping. Hence, 
by [He1, Theorem 0.5.7], 
there exists a long exact sequence

\begin{equation}
\label{chlong}
 \dots \rightarrow
H_{n+1}( {\cal L}_{\sim} ) \stackrel {\zeta_{n+1}} {\rightarrow}  H_n( {\cal K}_{\sim} ) 
 \stackrel {H_n({\psi}_{\sim})} {\rightarrow}  
 H_n({\cal P}_{\sim} ) \stackrel {H_n({\phi}_{\sim})} {\rightarrow} 
 H_n({\cal L}_{\sim}) \stackrel {\zeta_{n}} {\rightarrow}
 H_{n-1}({\cal K}_{\sim} ) \rightarrow \dots 
\end{equation}
of homology groups,
where $\zeta_n $ (the connecting morphism), $H_n({\phi}_{\sim})$ and $H_n({\psi}_{\sim})$ 
are continuous linear operators.

     We can also consider  a short sequence of dual cochain complexes
\[
  0 \rightarrow   {\cal L}_{\sim}^*  \stackrel {{\phi}_{\sim}^*} {\rightarrow}     
 {\cal P}_{\sim}^*  \stackrel {{\psi}_{\sim}^*} {\rightarrow}   {\cal K}_{\sim}^*
 \rightarrow   0
\]
in the category of Banach spaces and continuous linear operators.
This sequence is exact by [Jo1; Section 1] 
(see also, [He1; 0.5.2]). Hence, by [He1, Theorem 0.5.7], there exists a long exact sequence

\begin{equation}
\label{colong}
 \dots \rightarrow
H^{n-1}( {\cal K}_{\sim}^* ) \stackrel {\xi_n} {\rightarrow} H^n( {\cal L}_{\sim}^* ) 
 \stackrel {H^n({\phi}_{\sim}^*)} {\rightarrow}  
 H^n({\cal P}_{\sim}^* ) \stackrel {H^n({\psi}_{\sim}^*)} {\rightarrow}  
 H^n({\cal K}_{\sim}^*) \stackrel {\xi_{n+1}} {\rightarrow}
 H^{n+1}({\cal L}_{\sim}^* ) \dots 
\end{equation}
of cohomology groups, where $\xi_n$ (the connecting morphism),  $H^n({\phi}_{\sim}^*)$ and 
$H^n({\psi}_{\sim}^*)$ 
are continuous linear operators.

 (I) Since  $H^k({\psi}_{\sim}^*)$ is an isomorphism
whenever $n+2 \ge k \ge n-1$, by Lemma 3.1, $ H^k( {\cal L}_{\sim}^* ) = \{ 0 \}$ 
whenever  $n +2 \ge k \ge n$. Thus, by  [Jo1; Section 1],
 $ H_k( {\cal L}_{\sim} ) = \{ 0 \}$ 
for  $k =n $ and $k = n+1,$ and so  $H_n({\psi}_{\sim})$ is an 
 isomorphism. Hence, by  Lemma 7.1.32 [BS], 
$H_n({\psi}_{\sim})$ is a 
topological isomorphism. 

  (II)  Since  $H_k({\psi}_{\sim})$ is an isomorphism whenever 
$n+1 \ge k \ge n-2$, by Lemma 3.1, $ H_k( {\cal L}_{\sim}) = \{ 0 \}$ 
whenever  $n +1 \ge k \ge n-1$. Thus, by  [Jo1; Section 1],
 $ H^k( {\cal L}_{\sim}^* ) = \{ 0 \}$ 
for  $k =n $ and $k = n+1,$ and so  $H^n({\psi}_{\sim}^*)$ is an 
 isomorphism.  Hence, by Lemma 7.1.32 [BS], 
 $H^n({\psi}_{\sim}^*)$ is a 
topological isomorphism.

   (III) Recall that a complex 
$ {\cal L}_{\sim}$ is exact if and only if 
its dual complex $ {\cal L}_{\sim}^*$ is exact  [Jo1; Section 1]
 (see also [He1; 0.5.2]).
 Thus the homology groups
 $H_n( {\cal L}_{\sim} ) $ vanish if and only if the cohomology groups
$H^n({\cal L}_{\sim}^* ) $ vanish.

  Since the sequences (\ref{chlong}) and (\ref{colong}) are exact,
 the result follows from Lemma 3.1 
and Theorem 0.5.10 [He1]. 
\qed
 
 The following lemma is widely known.

{\sc 3.3 Lemma.}~~ {\it Let
\begin{equation}
\label{shexse}
 0 \rightarrow Y \stackrel {i} {\rightarrow}    Z  \stackrel {j} {\rightarrow}  W
 \rightarrow  0
\end{equation}
be a weakly admissible sequence of Banach spaces and continuous operators. 
 Then, for every Banach space $X$,  the sequence
\[
 0 \rightarrow  X \hat{\otimes} Y \stackrel {id_X \otimes i} {\rightarrow}    X \hat{\otimes} Z  
 \stackrel {id_X \otimes j} {\rightarrow}    X \hat{\otimes} W
 \rightarrow  0
\]
is  weakly admissible.}

 {\it Proof.} Since the sequence  (\ref{shexse})  is weakly admissible, there are 
continuous linear operators
$\alpha: Y^* \rightarrow   Z^*$ and $\beta: Z^* \rightarrow  
W^*$ such that $i^* \circ  \alpha = id_{Y^*}$, $\beta \circ j^* = id_{W^*}$ and
$ \alpha \circ i^*  +  j^* \circ \beta  = id_{Z^*}$. 

 Further we will use the well-known isomorphism (see, for example, 
[He1; Theorem 2.2.17])
\[
{\cal B}(X, Y^*)  \rightarrow (X \hat{\otimes} Y)^*: \phi \mapsto \Phi_\phi
\]
where $\Phi_\phi(x \otimes y) = [\phi(x)](y); x \in X, y \in Y$ and
\[
 (X \hat{\otimes} Y)^* \rightarrow {\cal B}(X, Y^*): f \mapsto \phi_f
\]
where $[\phi_ f(x)](y) = f(x \otimes y); x \in X, y \in Y.$

  We can see that, for  $ f \in (X \hat{\otimes} Y)^*  $, we have 
$\phi_ f \in {\cal B}(X, Y^*)  $,~~
$\alpha \circ \phi_ f \in {\cal B}(X, Z^*) $ and $\Phi_{\alpha \circ \phi_ f } \in 
(X \hat{\otimes} Z)^*  $. We define $ \gamma(f)
 = \Phi_{\alpha \circ \phi_ f}$ for  $ f \in (X \hat{\otimes} Y)^*  $ and $ \eta(g)
 = \Phi_{\beta \circ \phi_ g}$ for $ g \in (X \hat{\otimes} Z)^*  $.

 We can check that $(id_X \otimes i)^* \circ  \gamma = id_{(X \hat{\otimes} Y)^*}$,~~ 
$\eta \circ (id_X \otimes j)^* = id_{(X \hat{\otimes} W)^*}$ and
$ \gamma \circ (id_X \otimes i)^*  +  (id_X \otimes j)^* \circ \eta  = id_{(X \hat{\otimes} Z)^*}$. 
The result follows.
\qed

{\sc 3.4 Lemma.}~~ {\it Let 
\[
 0 \rightarrow Y \stackrel {i} {\rightarrow}    Z  \stackrel {j} {\rightarrow} W
 \rightarrow  0
\]
be a weakly admissible sequence of Banach spaces and continuous operators. 
Then, for every $n$,  the following 
 conditions are satisfied 

  {\rm (i)}  the operator 
\[
i^{\hat{\otimes} n} = i \hat{\otimes } \dots \hat{\otimes} i: 
Y \hat{\otimes} \dots \hat{\otimes} Y \rightarrow Z \hat{\otimes} \dots \hat{\otimes} Z
\]
is injective and
${\rm Im} (i^{\hat{\otimes} n})$ is closed;

 {\rm (ii)}  ${\rm Ker}~ {j^{\hat{\otimes} n}} =
Z^{\hat{\otimes} (n-1)} \hat{\otimes} i( Y ) +
 Z^{\hat{\otimes} (n-2)} \hat{\otimes} i( Y ) \hat{\otimes} Z + \dots +
i( Y ) \hat{\otimes}  Z^{\hat{\otimes} (n-1)}
$ and so it  is the closure of linear span 
$\{z_1 \otimes z_2 \otimes \dots \otimes z_n $: {\rm at least one} $z _i$ {\rm belongs} 
to $i(Y) \}.$

 {\rm (iii)}  the complex
\[
 0 \rightarrow  Y^{\hat{\otimes} n} \stackrel {i^{\hat{\otimes} n}} {\rightarrow} 
 Z^{\hat{\otimes} n}  \stackrel {j^{\hat{\otimes} n}} {\rightarrow} 
 (W)^{\hat{\otimes} n}
 \rightarrow  0
\]
is exact at $ Y^{\hat{\otimes} n}$  and $(W)^{\hat{\otimes} n}$;
 }

 {\it Proof.}  (i) We can see that
\[
 i^{\hat{\otimes} n} = 
 \kappa_n \circ (id_{Z^{\hat{\otimes} (n-1)}} \otimes i)  \circ
\kappa_{n-1} \circ \dots  
\circ (id_{Z \hat{\otimes} Y^{\hat{\otimes} (n-2)}} \otimes i) 
 \circ \kappa_1 \circ (id_{ Y^{\hat{\otimes} (n-1)}} \otimes i) ,
\]
where
$
\kappa_i: Z^{\hat{\otimes} (i-1)} \hat{\otimes} Y^{\hat{\otimes} (n-i)}
 \hat{\otimes} Z  \cong  
 Z^{\hat{\otimes} i} \hat{\otimes} Y^{\hat{\otimes} (n-i)}, \; i= 1,
\dots, n$
are given by $\kappa_i(x_1 \otimes x_2 \otimes \dots \otimes x_n) = 
x_n \otimes x_1 \otimes \dots \otimes x_{n-1}.$

  Thus 
\[
 (i^{\hat{\otimes} n})^* = (id_{ Y^{\hat{\otimes} (n-1)}} \otimes i)^* \circ (\kappa_1)^*
\circ (id_{Z \otimes Y^{\hat{\otimes} (n-2)}} \otimes i)^* \circ \dots \circ 
(\kappa_{n-1})^* \circ (id_{Z^{\hat{\otimes} (n-1)}} \otimes i)^*  \circ
 (\kappa_n)^*.
\]
By Lemma 3.3, $(i^{\hat{\otimes} n})^*$ is surjective as the composition of surjective 
operators. Therefore, by 
 [Ed; 8.6.15], $i^{\hat{\otimes} n}$ is injective and ${\rm Im} (i^{\hat{\otimes} n})$
 is closed in $Z^{\hat{\otimes} n} $.

 (ii) We can see that
\[
 j^{\hat{\otimes} n} = \psi_1  \circ \psi_2  \circ \dots  \circ \psi_n
\]
where $ \psi_1 =  j \otimes  id_{W^{\hat{\otimes} (n-1)}}$, 
 $ \psi_i = id_{Z^{\hat{\otimes} (i-1)}} \otimes j \otimes  id_{W^{\hat{\otimes} (n-i)}}$
for $ 1 < i < n $  and $ \psi_n = id_{Z^{\hat{\otimes} (n-1)}} \otimes j .$
By Lemma 3.3,  we have ${\rm Ker}~ \psi_i = {\rm Im}~
 id_{Z^{\hat{\otimes} (i-1)}} \otimes i \otimes  
id_{W^{\hat{\otimes} (n-i)}} =
Z^{\hat{\otimes} (i-1)} \hat{\otimes}i( Y ) \hat{\otimes} W^{\hat{\otimes} (n-i)}$ and, 
for $ F_i = Z^{\hat{\otimes} (i-1)} \hat{\otimes} i( Y )  \hat{\otimes} Z^{\hat{\otimes} (n-i)}$,
we obtain
\[
(\psi_{i+1}  \circ \dots  \circ \psi_{n-1}  \circ \psi_n) (F_i) 
= {\rm Ker}~ \psi_i.
\]
Thus we have
\[
{\rm Ker}~ \psi_n =  Z^{\hat{\otimes} (n-1)}  \hat{\otimes} i( Y );
\]
\[
{\rm Ker}~ (\psi_{n-1} \circ \psi_n) = 
\psi_{n}^{-1} ({\rm Ker}~\psi_{n-1}) =
\]
\[
{\rm Ker}~ \psi_n + F_{n-1}  = 
 Z^{\hat{\otimes} (n-1)}  \hat{\otimes} i( Y )  + 
 Z^{\hat{\otimes} (n-2)}  \hat{\otimes} i( Y ) \hat{\otimes} Z ;
\]
\[
{\rm Ker}~ (\psi_{n-2} \circ \psi_{n-1} \circ \psi_n) = 
(\psi_{n-1} \circ \psi_{n})^{-1} ({\rm Ker}~\psi_{n-2}) =
\]
\[
{\rm Ker}~ (\psi_{n-1} \circ \psi_n)  +  F_{n-2}   = 
\]
\[
 Z^{\hat{\otimes} (n-1)}  \hat{\otimes}  i( Y )  + 
 Z^{\hat{\otimes} (n-2)}  \hat{\otimes} i( Y ) \hat{\otimes} Z +
Z^{\hat{\otimes} (n-3)}  \hat{\otimes} i( Y ) \hat{\otimes} Z^{\hat{\otimes}
2};
\]
\[
\dots \dots 
\]
\[
{\rm Ker}~ (\psi_1  \circ \dots \psi_{n-1} \circ \psi_n) = 
(\psi_2  \circ \dots \psi_{n-1} \circ \psi_n)^{-1} ({\rm Ker}~(\psi_{1}) =
\]
\[
{\rm Ker}~ (\psi_2  \circ \dots \psi_{n-1} \circ \psi_n) + F_{1} =
\]
\[
Z^{\hat{\otimes} (n-1)} \hat{\otimes} i( Y ) +
 Z^{\hat{\otimes} (n-2)} \hat{\otimes} i( Y ) \hat{\otimes} Z + \dots +
i( Y ) \hat{\otimes}  Z^{\hat{\otimes} (n-1)}.
\]

 (iii) It follows from (i) and the fact that the projective tensor product of surjective operators is
also surjective.
\qed

{\sc 3.5 Theorem.}~~ {\it Let
\begin{equation}
\label{waeb}
 0 \rightarrow B \stackrel {i} {\rightarrow}    A  \stackrel {j} {\rightarrow} D
 \rightarrow  0
\end{equation}
be a weakly admissible extension of Banach algebras. Then 

 {\rm (I)} there exists an associated long exact
sequence of Banach simplicial homology groups 
\begin{equation}
\label{lessh}
\dots  \rightarrow {\cal H}_n(B) \rightarrow
{\cal H}_n(A) \rightarrow  {\cal H}_n(D) \rightarrow
 {\cal H}_{n-1}(B) \rightarrow \dots \rightarrow
{\cal H}_0(D) \rightarrow 0
\end{equation}
if and only if there exists an associated long exact
sequence of Banach simplicial cohomology groups 
\begin{equation}
\label{lessc}
0 \rightarrow {\cal H}^0(D)   \rightarrow
\dots \rightarrow
{\cal H}^{n-1}(B) \rightarrow  {\cal H}^n(D) \rightarrow
 {\cal H}^n(A) \rightarrow {\cal H}^n(B) \rightarrow \dots;
\end{equation}

{\rm (II)}  there exists an associated long exact
sequence of Banach Bar homology groups 
\begin{equation}
\label{lesBh}
\dots \rightarrow  {\cal HR}_n(B) \rightarrow
{\cal HR}_n(A) \rightarrow  {\cal HR}_n(D) \rightarrow
 {\cal HR}_{n-1}(B) \rightarrow \dots \rightarrow
{\cal HR}_0(D) \rightarrow 0
\end{equation}
if and only if  there exists an associated long exact
sequence of Banach Bar cohomology groups 
\begin{equation}
\label{lesBc}
 0 \rightarrow {\cal HR}^0(D)  \rightarrow
 \dots \rightarrow
{\cal HR}^{n-1}(B) \rightarrow  {\cal HR}^n(D) \rightarrow
 {\cal HR}^n(A) \rightarrow {\cal HR}^n(B) \rightarrow \dots .
\end{equation}
}

{\it Proof.} (I) By Lemma 3.4 (iii), for  $n = 0, 1, \dots,$  the extension  (\ref{waeb})
induces a complex
of Banach spaces and continuous operators
\[
  0 \rightarrow  C_n(B)  \stackrel {i^{\hat{\otimes} (n+1)}} {\rightarrow} 
 C_n(A)  \stackrel {j^{\hat{\otimes} (n+1)}} {\rightarrow}  C_n(D) 
 \rightarrow   0
\]
which is exact at $ C_n(D) $ and $  C_n(B) $.
It is routine to check that the diagram
\[
\begin{array}{ccccccccc}
~ & ~ &  0 & ~ &   0 & ~ &  0 & ~ & ~
\\
~ & ~ &  \uparrow & ~ &   \uparrow & ~ &  \uparrow & ~ & ~
\\
    0  & \rightarrow & C_0(B) & \stackrel {i} {\rightarrow} &
 C_0(A)  & \stackrel {j} {\rightarrow} & C_0(D) &
 \rightarrow &  0
\\
~ & ~ &  \uparrow d_0 & ~ &   \uparrow d_0 & ~ &  
\uparrow d_0 & ~ & ~
\\
    0  & \rightarrow & C_1(B) & \stackrel {i \hat{\otimes} i} {\rightarrow} &
 C_1(A)  & \stackrel {j \hat{\otimes} j} {\rightarrow} & C_1(D) &
 \rightarrow &  0
\\
~ & ~ &  \uparrow d_1 & ~ &   \uparrow d_1 & ~ &  
\uparrow d_1 & ~ & ~
\\
\dots & \dots &  \dots  & \dots &   \dots & \dots &  \dots & \dots & \dots
\\
~ & ~ &  \uparrow & ~ &   \uparrow & ~ &  \uparrow & ~ & ~
\\
    0  & \rightarrow & C_n(B) & \stackrel {i \hat{\otimes} \dots \hat{\otimes} i} {\rightarrow} &
 C_n(A)  & \stackrel {j \hat{\otimes} \dots \hat{\otimes} j} {\rightarrow} & C_n(D) &
 \rightarrow &  0
\\
~ & ~ &  \uparrow d_n & ~ &   \uparrow d_n & ~ &  
\uparrow d_n  & ~ & ~
\\
    0  & \rightarrow & C_{n+1}(B) & \stackrel {i \hat{\otimes} \dots \hat{\otimes} i} {\rightarrow} &
 C_{n+1}(A)  & \stackrel {j \hat{\otimes} \dots \hat{\otimes} j} {\rightarrow} & C_{n+1}(D) &
 \rightarrow &  0
\\
~ & ~ &  \uparrow & ~ &   \uparrow & ~ &  \uparrow & ~ & ~
\\
\dots & \dots &  \dots  & \dots &   \dots & \dots &  \dots & \dots & \dots
\end{array}
\]
is commutative. The same is true for cochains.
Thus the extension (\ref{waeb}) induces a
short exact sequence of chain complexes
\[
  0 \rightarrow {\cal C}_{\sim}(A,D)  \rightarrow     
 {\cal C}_{\sim}(A)  \rightarrow   {\cal C}_{\sim}(D)
 \rightarrow   0
\]
and
a short exact sequence of cochain complexes
\[
  0 \rightarrow   {\cal C}^{\sim}(D)  \rightarrow     
 {\cal C}^{\sim}(A)  \rightarrow     {\cal C}^{\sim}(A; D)
 \rightarrow   0
\]
in the category of Banach spaces and continuous operators, where
$ {\cal C}_{\sim}(A,D) $ is the subcomplex ${\rm Ker }(j \hat{\otimes} \dots \hat{\otimes} j)$ of
$ {\cal C}_{\sim} (A)$; 
$ {\cal C}^{\sim}(A; D) $ is the subcomplex 
$[{\rm Ker }(j \hat{\otimes} \dots \hat{\otimes} j)]^*$ of
$ {\cal C}^{\sim}(A) $. Hence, by [He1, Chapter 0, 
Section 5.5], there exist a long exact sequence of homology groups
\[
 \dots \rightarrow
{\cal H}_{n+1}(D)  \rightarrow  H_n({\cal C}_{\sim}(A; D) )
 \rightarrow
 {\cal H}_n(A) \rightarrow {\cal H}_n(D) 
 \rightarrow \dots 
\]
and a
 long exact sequence  of cohomology groups
\[
 \dots \rightarrow
H^{n-1}({\cal C}^{\sim}(A; D) ) \rightarrow  {\cal H}^n(D) \rightarrow
 {\cal H}^n(A) \rightarrow  H^n({\cal C}^{\sim}(A; D) ) \rightarrow \dots.
\]

By Lemma 3.4 (i),
\[
i^{\hat{\otimes} (n+1)}: 
B^{\hat{\otimes} (n+1)}  \rightarrow  {\rm Ker}~ j^{\hat{\otimes} (n+1)}
\]
is injective and
${\rm Im} (i^{\hat{\otimes} (n+1)})$ is closed.
Therefore, by Proposition 3.2, for the two complexes ${\cal C}_{\sim}(B)$ and 
${\cal C}_{\sim}(A; D)$
and the continuous morphism  $\{i^{\hat{\otimes} (n+1)}\}$ of complexes, we 
have  
 $ H_n({\cal C}_{\sim}(A; D) ) \cong {\cal H}_n(B) $ for  $n = 0, 1, \dots$
 if and only if 
$ H^n({\cal C}^{\sim}(A; D) ) \cong {\cal H}^n(B) $ for $n = 0, 1, \dots$.
 The result follows from the five lemma [Ma, Lemma 1.3.3].

 (II) follows by the same arguments as in the case (I).
\qed

 The following lemma is widely known.

{\sc 3.6 Lemma.}~~ {\it Let
\begin{equation}
\label{exten}
 0 \rightarrow B \stackrel {i} {\rightarrow}    A  \stackrel {j} {\rightarrow} D
 \rightarrow  0
\end{equation}
 be an extension of Banach algebras. Suppose that $B$ has a left or right bounded approximate
identity $(e_{\alpha})_{ \alpha \in \Lambda}$. Then } (\ref{exten}) {\it is a weakly admissible.}

 {\it Proof.}  Consider the Fr${\rm \acute{e}}$chet filter $F$ on $\Lambda$, with 
 base $ \{ Q_{\lambda}: 
\lambda \in \Lambda \}$, where 
$Q_{\lambda} = \{ \alpha \in \Lambda : \alpha \ge \lambda \}.$ Thus
\[
 F = 
\{ E \subset \Lambda : \;  {\rm there \; is \; a} \; \lambda  \in \Lambda \; {\rm such \; that } \; 
 Q_{\lambda} \subset E \}.
\]
Let $U$ be an ultrafilter on $ \Lambda $ which refines $F$. One can find information on 
 filters in [Bo]. 
   
  Suppose that $(e_{\alpha})_{ \alpha \in \Lambda}$ is a  right bounded approximate
identity. For $f \in B^*$ we define $g_f \in A^*$ by
\[
g_f(a) =\lim_{ \alpha \rightarrow U} f(i^{-1}( a i(e_{\alpha}))) \; 
{\rm for \; all} \; a \in A.
\]
It is easy to check that $g_f$ is a bounded linear functional, the operator 
\[
L : B^* \rightarrow A^*: f \mapsto g_f
\]
is a bounded linear operator and  $  i^* \circ L = id _{B^*}$.
\qed

  Let 
\[
 0 \rightarrow B \stackrel {i} {\rightarrow}    A  \stackrel {j} {\rightarrow} D
 \rightarrow  0
\]
be an extension of Banach algebras. 
 Note that
$(id_{A^{\hat{\otimes} (n+1)}} - t_{\underline{n}})
 \circ i^{\hat{\otimes} (n+1)} =
  i^{\hat{\otimes} (n+1)} \circ (id_{B^{\hat{\otimes} (n+1)}}
 - t_{\underline{n}})$ and 
$(id_{D^{\hat{\otimes} (n+1)}} - t_{\underline{n}}) 
\circ j^{\hat{\otimes} (n+1)} =
  j^{\hat{\otimes} (n+1)} \circ (id_{A^{\hat{\otimes} (n+1)}} 
- t_{\underline{n}})$.
 Thus there is a complex
\[
 0 \rightarrow  CC_n(B) \stackrel {\widetilde{i^{\hat{\otimes} (n+1)}}} {\rightarrow} 
  CC_n(A)  \stackrel {\widetilde{j^{\hat{\otimes} (n+1)}}} {\rightarrow} 
  CC_n(D)
 \rightarrow  0,
\]
where $\widetilde{i^{\hat{\otimes} (n+1)}}$ and $\widetilde{j^{\hat{\otimes} (n+1)}}$ are induced by 
 $i^{\hat{\otimes} (n+1)}$ and $j^{\hat{\otimes} (n+1)}$ respectively.

{\sc 3.7 Lemma.}~~ {\it Let 
\[
 0 \rightarrow B \stackrel {i} {\rightarrow}    A  \stackrel {j} {\rightarrow} D
 \rightarrow  0
\]
be an extension of Banach algebras. 
Suppose $B$ has a left or right bounded approximate identity 
$(e_{\alpha})_{ \alpha \in \Lambda}$.
Then, for every $n$,  the complex
\[
 0 \rightarrow  CC_n(B) \stackrel {\widetilde{i^{\hat{\otimes} (n+1)}}} {\rightarrow} 
  CC_n(A)  \stackrel {\widetilde{j^{\hat{\otimes} (n+1)}}} {\rightarrow} 
  CC_n(D)
 \rightarrow  0,
\]
is exact at $  CC_n(B)$  and $ CC_n(D)$;  $\widetilde{i^{\hat{\otimes} (n+1)}}$ is injective and
${\rm Im} (\widetilde{i^{\hat{\otimes} (n+1)}})$ is closed; 
 and  ${\rm Ker}~ {\widetilde{j^{\hat{\otimes} (n+1)}}} = \{u + 
{\rm Im}~(id_{A^{\hat{\otimes} (n+1)}} -   t_{\underline{n}});$
where $ u \in {\rm Ker}~j^{\hat{\otimes} (n+1)} \}$
and so it 
is the closure of the linear span of 
$\{a_0 \otimes a_1 \otimes \dots \otimes a_n + 
{\rm Im}~(id_{A^{\hat{\otimes} (n+1)}} -   t_{\underline{n}}) :$ 
{\rm  at least one} $ a_i$ {\rm belongs to} $i(B) \}.$
}

 {\it Proof.}  Since $ j^{\hat{\otimes} (n+1)}$ is surjective and 
 $(id_{D^{\hat{\otimes} (n+1)}} - t_{\underline{n}}) \circ j^{\hat{\otimes} (n+1)} =
  j^{\hat{\otimes} (n+1)} \circ (id_{A^{\hat{\otimes} (n+1)}} - t_{\underline{n}})$,
the mapping 
$ \widetilde{j^{\hat{\otimes} (n+1)}}: CC_n(A) \rightarrow CC_n(D)$ is surjective too.

  Suppose that $(e_{\alpha})_{ \alpha \in \Lambda}$ is a  right bounded approximate
identity. Let us check that $(\widetilde{i^{\hat{\otimes} (n+1)}})^*$ is surjective.
For
\[
 f \in CC_n(B)^* = \{ f \in (B^{\hat{\otimes} (n+1)})^* :
 f|_{{\rm Im}~(id_{B^{\hat{\otimes} (n+1)}} -
t_{\underline{n}})}
 = 0 \},
\]
we define $g_f \in (A^{\hat{\otimes} (n+1)})^*$ by
\[
g_f(a_0 \otimes a_1 \otimes \dots \otimes a_n) =
\lim_{ \alpha \rightarrow U} 
 f(i^{-1}( a_0 i(e_{\alpha})) \otimes \dots \otimes i^{-1}( a_n i(e_{\alpha}))),
\]
where $a_i \in A;~~ i= 0, \dots, n$ and the ultrafilter $U$ is 
defined as in Lemma 3.6.
It is easy to check that $g_f$ is a bounded linear functional,
$g_f|_{ {\rm Im}~(id_{A^{\hat{\otimes} (n+1)}} 
- t_{\underline{n}})} = 0$ and
 $f = (i^{\hat{\otimes} (n+1)})^* g_f.$
Thus, by [Ed; 8.6.15], 
$\widetilde{i^{\hat{\otimes} (n+1)}}: CC_n(B) \rightarrow CC_n(A)$ is  injective  and its image is
closed.

  Now let us show that there exists a continuous operator 
$\alpha_n: C_n(A)^* \rightarrow C_n(D)^*$  such that 
$\alpha_n \circ (j^{\hat{\otimes} (n+1)})^*  = id_{C_n(D)^*}$ and 
$\alpha_n(CC_n(A)^*) \subset CC_n(D)^*$. 
For
$ g \in C_n(A)^*$, 
we define $f_g \in C_n(D)^*$, by
\[
f_g(d_0 \otimes d_1 \otimes \dots \otimes d_n) =
\]
\[
g(a_0 \otimes a_1 \otimes \dots \otimes a_n) -
\lim_{ \alpha \rightarrow U} \sum_{k=0}^{n} 
 g( a_0 \otimes \dots   \otimes a_{k-1} \otimes a_k i(e_{\alpha})\otimes
 a_{k+1} \otimes \dots \otimes a_n) +
\]
\[
\lim_{ \alpha \rightarrow U} \sum_{0 \le k< t \le n}
 g( a_0 \otimes \dots   \otimes a_{k-1} \otimes a_k i(e_{\alpha})\otimes
 a_{k+1} \otimes \dots \otimes a_t i(e_{\alpha})\otimes a_{t+1}\otimes \dots 
\otimes a_n) - \dots +
\]
\[
(-1)^{n+1} \lim_{ \alpha \rightarrow U} 
 g( a_0 i(e_{\alpha}) \otimes  a_1 i(e_{\alpha}) \otimes \dots 
\otimes a_n i(e_{\alpha})),
\]
where $d_i = j(a_i) \in D;~~ i= 0, \dots, n$.
It is easy to check that the operator
\[
\alpha_n: C_n(A)^* \rightarrow C_n(D)^*: g \mapsto f_g
\]
is a bounded linear operator,~~
$\alpha_n \circ (j^{\hat{\otimes} (n+1)})^* = id_{C_n(D)^*}$ and 
$\alpha_n(CC_n(A)^*) \subset CC_n(D)^*$. 

  Hence we can see that the dual map of 
\[
\sigma_n: {\rm Ker}~j^{\hat{\otimes} (n+1)} \rightarrow 
{\rm Ker}~\widetilde{j^{\hat{\otimes} (n+1)}}: u \mapsto u +
{\rm Im}~(id_{A^{\hat{\otimes} (n+1)}} -   t_{\underline{n}})
\]
is injective and ${\rm Im}~\sigma_n$ is closed. Thus,
by [Ed; 8.6.15],  $\sigma_n$ is surjective, that is, 
 $({\rm Ker}~\widetilde{j^{\hat{\otimes} (n+1)}})^* ={\rm Im}~\sigma_n$.
It follows from Lemma 3.4(ii) that  the ${\rm Ker}~ {\widetilde{j^{\hat{\otimes} (n+1)}}}$ 
is the closure of the linear span of
$\{a_0 \otimes a_1 \otimes \dots \otimes a_n + 
{\rm Im}~(id_{A^{\hat{\otimes} (n+1)}} -   t_{\underline{n}}) :$ 
{\rm  at least one} $a _i$ {\rm belongs to} $i(B) \}.$ 
\qed

{\sc 3.8 Theorem.}~~ {\it Let
\[
 0 \rightarrow B \stackrel {i} {\rightarrow}    A  \stackrel {j} {\rightarrow} D
 \rightarrow  0
\]
 be an extension of Banach algebras.
 Suppose $B$ has a left or right bounded approximate identity.
Then  there exists an associated long exact
sequence of Banach  cyclic homology groups 
\begin{equation}
\label{lesch}
\dots \rightarrow {\cal HC}_n(B) \rightarrow
{\cal HC}_n(A) \rightarrow  {\cal HC}_n(D) \rightarrow
 {\cal HC}_{n-1}(B) \rightarrow \dots \rightarrow
{\cal HC}_0(D) \rightarrow 0
\end{equation}
if and only if  there exists an associated long exact
sequence of Banach  cyclic cohomology groups 
\begin{equation}
\label{lescc}
0 \rightarrow {\cal HC}^0(D)  \rightarrow
 \dots \rightarrow
{\cal HC}^{n-1}(B) \rightarrow  {\cal HC}^n(D) \rightarrow
 {\cal HC}^n(A) \rightarrow {\cal HC}^n(B) \rightarrow \dots .
\end{equation}
}

{\it The proof} requires only minor modifications of that of 
Theorem 3.5 in view
of Lemmas 3.6 and 3.7.
\qed

\section{The existence of long exact sequences of cohomology groups 
of Banach algebras }

  As stated above, Wodzicki remarked 
that his result on the existence of long exact sequences of cyclic and simplicial
 homology groups can be extended to the continuous case.
 In Corollary 4 [Wo2] he considered extensions of Banach 
algebras
$ 0 \rightarrow B \stackrel {i} {\rightarrow}    A  \stackrel {j} {\rightarrow} D
 \rightarrow  0$
where $B$ has a left or right bounded approximate identity and stated 
 the existence of  long exact sequences 
associated with the extension under the explicit hypothesis of
the admissibility of the extension. 
Subsequently, he pointed out that 
this hypothesis is unnecessary.
Here we state his result in a new formulation.
   
{\sc 4.1 Theorem.}~~ {\it Let
\[
0 \rightarrow B \stackrel {i} {\rightarrow}    A  \stackrel {j} {\rightarrow} D
 \rightarrow  0 
\]
 be an extension of Banach algebras. 
Suppose $B$ has a left or right bounded approximate identity.
 Then  there exist  associated long exact
sequences  of Banach simplicial homology groups (\ref{lessh}),
 of Banach  cyclic homology groups  (\ref{lesch})
and  of Banach Bar homology groups  (\ref{lesBh}),
and so ${\cal HR}_n(A) = {\cal HR}_n(D)$ for all $n \ge 0$.
}

{\it Proof.} Let us introduce the
following notation. Let 
 $\underline{n} = (n_1, \dots, n_{l+1})$ be an
$(l+1)$-tuple of integers such that 
$n_1,n_{l+1} \ge 0$ 
and the others are $> 0$ and let  $\underline{k} = (k_1, \dots,
k_{l})$ be an
$l$-tuple of integers such that 
$k_1 \ge 0$ 
and the others are $> 0$.
 Put  $|\underline{n}|: = n_1+ \dots + n_{l+1}$ 
and $l(\underline{n}): = l$.
 In view of Lemmas 3.3, 3.4, 3.6 and 3.7, we can see that,
up to algebraic isomorphism,
\[
C_m(A) = \bigoplus_{\underline{n}, \underline{k}: 
|\underline{n}| +|\underline{k}| = m+1, \; 
l(\underline{n})\ge 1, \; l(\underline{k}) = l(\underline{n})-1  } 
D^{\hat{\otimes}n_1} \hat{\otimes} B^{\hat{\otimes}k_1} \hat{\otimes}
\dots \hat{\otimes} D^{\hat{\otimes}n_{l+1}}, 
\] 
\[
{\rm Ker}~j^{\hat{\otimes}(m+1)} = \bigoplus_{\underline{n},\;
 \underline{k}: k_1 > 0, \;
|\underline{n}| +|\underline{k}| = m+1, \; 
l(\underline{n})\ge 1, \; l(\underline{k}) = l(\underline{n})-1  } 
D^{\hat{\otimes}n_1} \hat{\otimes} B^{\hat{\otimes}k_1} \hat{\otimes}
\dots \hat{\otimes} D^{\hat{\otimes}n_{l+1}}
\] 
 and  ${\rm Ker}~ {\widetilde{j^{\hat{\otimes} (m+1)}}} = \{u + 
{\rm Im}~(id_{A^{\hat{\otimes} (m+1)}} -   t_{\underline{m}});$
where $ u \in {\rm Ker}~j^{\hat{\otimes} (m+1)} \}$.
Thus  the  proof is the same as in the algebraic case
in Theorem 3.1 [Wo3] and Theorem 3  [Wo2]
 with filtrations in which algebraic tensor products are replaced by the
projective tensor products.
\qed

{\sc 4.2 Theorem.}~~{\it Let
\[
 0 \rightarrow B \stackrel {i} {\rightarrow}    A  \stackrel {j} {\rightarrow} D
 \rightarrow  0 
\] 
 be an extension of Banach algebras. 
Suppose $B$  has a left or right bounded approximate identity.
 Then  there exist associated long exact
sequences of Banach simplicial cohomology groups (\ref{lessc}),
of Banach  cyclic cohomology groups  (\ref{lescc})
and of Banach Bar cohomology groups (\ref{lesBc}),
and so $  {\cal HR}^n(D) = {\cal HR}^n(A)$ for all $ n \ge 0$. 
}

{\it Proof.} The existence of the long exact sequences follows 
from  Theorems 4.1, 3.5
and 3.8 and Lemma 3.6.

   By Theorem 16 [He2],  $  {\cal HR}^n(D) = \{ 0\}$ for 
every Banach algebra $A$
with left or right bounded approximate identity and for all $ n \ge 0$. Hence
the final statement follows.
\qed

     Since $C^*$-algebras have bounded approximate identities, 
Theorem 4.2 applies
whenever $B$ is a $C^*-$algebra. By  virtue of the main result of [Dix], 
the Banach algebra 
${\cal K}(E)$ of compact operators on a Banach space  $E$
with the bounded compact approximation property 
has a bounded  left approximate identity. For more examples of Banach algebras with a
left or right  bounded approximate identity see, for example  [Pa; Section
5.1].

{\sc 4.3 Remark.}~~ Let us consider the extension 
of Banach algebras
$ 0 \rightarrow B \stackrel {i} {\rightarrow}   A  \stackrel {j} {\rightarrow} 
D \rightarrow  0,$
when $B$ has a left or right bounded approximate identity. The existence of 
 the Connes-Tsygan sequence for one 
of Banach algebras $A$ or $D$ implies the existence of 
the same exact sequence for the other.
This follows from Theorems 15, 16 [He2] and Theorem 4.2.

{\sc 4.4 Proposition.}~~ {\it Let
\[
 0 \rightarrow B \stackrel {i} {\rightarrow}   A  \stackrel {j} {\rightarrow} 
D \rightarrow  0,
\]
be  an extension of Banach algebras. 
 Then

  ${\rm (i)}$ if $B$ is an amenable Banach algebra, then 
${\cal H}^n(A) = {\cal H}^n(D)  \; {\it
for \; all \;} n \ge 2 $
and there exist exact sequences
\[
0 \rightarrow D^{tr}   \rightarrow A^{tr}  \rightarrow
 B^{tr} 
 \rightarrow {\cal H}^{1}(D) \rightarrow  {\cal H}^{1}(A) 
\rightarrow 0 
\]
and
\[
  0 \rightarrow
{\cal HC}^{2n}(D) \rightarrow  {\cal HC}^{2n}(A) \rightarrow
  B^{tr}
  \rightarrow {\cal HC}^{2n+1}(D) \rightarrow  {\cal HC}^{2n+1}(A)
 \rightarrow 0 
\]
for every $n \ge 0$;

 ${\rm (ii)}$ if $B$ is a $C^*$-algebra without non-zero bounded 
traces then
\[
{\cal H}^n(A) = {\cal H}^n(D) \; {\it and}  \; 
{\cal HC}^n(A) = {\cal HC}^n(D) \; {\it
for \; all \;} n \ge 0.
\]
 } 

{\it Proof.}~~(i) By assumption $B$ is amenable, and so $B$ has a bounded approximate
identity. Hence, by Theorem 4.2, there 
 exist associated long exact sequences of Banach simplicial and 
cyclic cohomology groups (\ref{lessc}) and (\ref{lescc}).
Recall that, for an amenable Banach algebra $B$, ${\cal H}^n(B) =
\{ 0 \} $ for all $n \ge 1$ and, by  Theorem 25  [He2],
${\cal HC}^n(B)$ is equal
to  $B^{tr}$ for even $n$ and to $0$ for odd $n$.
The result follows in view of Lemma 7.1.32 [BS].

 (ii) By Proposition 1.7.3 [Di], any $C^*$-algebra  has  
a bounded approximate identity and so Theorem 4.2 applies.
 By Theorem 4.1 and Corollary 3.3 [ChS], for every $C^*$-algebras without
non-zero bounded traces,
 the Banach simplicial and cyclic cohomology groups vanish for all $n \ge 0$.
\qed

  Note that some results of this proposition were proved by a different approach in [Ly].

{\sc 4.5 Remark.}~~ Let us consider the extension 
of Banach algebras
$ 0 \rightarrow B \stackrel {i} {\rightarrow}   A  \stackrel {j} {\rightarrow} 
D \rightarrow  0,$ 
where $B$ has a left or right bounded approximate identity. Then the
 condition that $D$  be 
 an amenable Banach algebra  implies, by the same arguments as in 
Theorem 4.4,
  the following:
$
{\cal H}^n(A) = {\cal H}^n(B)  \; {\rm
for \; all \;} n \ge 1
$
and there exist exact sequences
\[
0 \rightarrow D^{tr}   \rightarrow A^{tr}  \rightarrow
 B^{tr}  \rightarrow 0
\]
and
\[
  0 \rightarrow
{\cal HC}^{2n+1}(A) \rightarrow  {\cal HC}^{2n+1}(B) \rightarrow
  D^{tr}
  \rightarrow {\cal HC}^{2n+2}(A) \rightarrow  {\cal HC}^{2n+2}(B)
 \rightarrow 0
\]
for every $n \ge 0$.

  As in Theorem 4.4, the  condition that $D$ be a $C^*$-algebra without non-zero bounded traces 
implies
${\cal H}^n(A) = {\cal H}^n(B) \; {\rm and}  \; 
{\cal HC}^n(A) = {\cal HC}^n(B) \; {\rm
for \; all \;} n \ge 0.$

  {\sc 4.6 Examples.}~~ Some examples of $C^*$-algebras without
non-zero bounded traces are: (i)  The $C^*$-algebra 
${\cal K}(H)$ of compact operators on  an infinite-dimensional Hilbert space $H$;
see [An; Theorem 2]. We can also show that $ C(\Omega,{\cal K}(H))^{tr} = 0$,
where $\Omega$ is a compact space.  (ii)  Properly infinite von Neumann
algebras ${\cal U}$. By Proposition 2.2.4 [Sa] in ${\cal U}$ there exists 
a sequence $(p_m)$ of mutually orthogonal, equivalent
projections  with $p_m \sim e$. Thus, by Theorem 2.1 [Fa], each hermitian
element of  ${\cal U}$ is the sum of five commutators. Hence there are no 
non-zero traces on ${\cal U}$. This class includes the   $C^*$-algebra 
${\cal B}(H)$ of all bounded operators on  an infinite-dimensional Hilbert space $H$;
see also [Hal] for the statement
${\cal B}(H)^{tr} = 0$.

{\sc 4.7 Examples.}~~ Each nuclear $C^*$-algebra is amenable [Ha]. 
Some examples of amenable $C^*$-algebras are: 
(i) $GCR~C^*$-algebras; in particular, commutative
$C^*$-algebras and the $C^*$-algebra of compact operators ${\cal K}(H)$ 
on a Hilbert space  $H$; (ii) Uniformly hyperfinite algebras (UHF-algebras)
[Mu; Remark 6.2.4].

  The group algebra $L^1(G)$ of Haar integrable functions on an amenable 
locally compact group $G$  with  convolution product  is amenable too [Jo1]. 
For more examples of amenable Banach algebras 
see, for example [GJW].

 {\sc 4.8  Examples.}~~ In  [GJW] it is shown that  the Banach algebra 
 ${\cal K}(E)$ of compact operators on a
Banach space $E$ with  property $({\Bbb A})$ which was defined in [GJW]
 is amenable. Property 
$({\Bbb A})$ implies that ${\cal K}(E)$ contains a bounded sequence of projections of
unbounded finite rank, and from this it is easy to show (via embedding of matrix
algebras) that there is no non-zero bounded trace on  ${\cal K}(E)$.
Thus we can see from Theorem 4.4 that, for every extension 
of Banach algebras
\[
 0 \rightarrow {\cal K}(E) \stackrel {i} {\rightarrow}   A  \stackrel {j} {\rightarrow} 
D \rightarrow  0,
\]
we have 
\[
{\cal H}^n(A) = {\cal H}^n(D) \; {\rm and}  \; 
{\cal HC}^n(A) = {\cal HC}^n(D) \; {\rm
for \; all \;} n \ge 0.
\]
In particular, for  the Banach algebra 
$A = {\cal B}(E)$ of all bounded operators on a
Banach space $E$,
${\cal H}^n({\cal B}(E)) = 
{\cal H}^n({\cal B}(E)/{\cal K}(E))$  and 
${\cal HC}^n({\cal B}(E)) = 
{\cal HC}^n({\cal B}(E)/{\cal K}(E))$ for all $n \ge 0$.
 Several classes of Banach spaces have the property $({\Bbb A})$: 
$l_p; 1 < p < \infty; C(K)$,
where $K$ is a compact Hausdorff space; $L_p (\Omega, \mu); 1 < p < \infty$
 (for details and more examples see [Jo1] and 
[GJW]).

\end{document}